\documentclass[a4paper,11pt]{amsart}
\usepackage{mathtext}         %

\usepackage{latexsym}
\usepackage{amsmath,amssymb,amsthm,amsopn,amscd}
\usepackage[dvips]{graphicx}
\usepackage{array,fullpage,wrapfig}

\usepackage[matrix,arrow,curve]{xy}
\usepackage[active]{srcltx}

\newcommand{\picref}[1]{\ref{#1}}

\newtheorem{theorem}[equation]{Theorem}
\newtheorem*{theorem*}{Theorem}
\newtheorem{proposition}[equation]{Proposition}
\newtheorem*{proposition*}{Proposition}

\newtheorem{lemma}[equation]{Lemma}
\newtheorem{corollary}[equation]{Corollary}
\newtheorem*{corollary*}{Corollary}

\theoremstyle{definition}

\newtheorem{definition}[equation]{Definition}

\DeclareMathOperator{\st}{st}

\theoremstyle{remark}
\newtheorem{remark}[equation]{Remark}

\begin{document}  

%commands
\newcommand{\ldot}{{\:\raisebox{1.5pt}{\selectfont\text{\circle*{1.5}}}}}
\newcommand{\udot}{{\:\raisebox{4pt}{\selectfont\text{\circle*{1.5}}}}}
\def\pt{\udot}

\def\kk{\Bbbk}

\let\le\leqslant
\let\ge\geqslant
\let\leq\leqslant
\let\geq\geqslant
\let\cong\simeq
\def\OP{\mathop\oplus\limits}
\def\OT{\mathop\otimes\limits}
\def\Sum{\mathop\sum\limits}
\newcommand{\bydef}{\stackrel{{\rm def}}{=}}
\newcommand{\Slut}{$\quad\Box$}
\newcommand{\Skip}{\smallskip\flushleft}
\newcommand{\Dim}{\textrm{dim}}
\newcommand{\Deg}{\textrm{deg}}
\newcommand{\Det}{\textrm{det}}
\newcommand{\al}{\alpha}
\newcommand{\be}{\beta}
\newcommand{\la}{\lambda}
\newcommand{\Tak}{\hat{\alpha}_{i}}
\newcommand {\bC} {\mathbb {C}}
\newcommand {\bR} {\mathbb R}
\newcommand {\bZ} {\mathbb Z}
\newcommand {\D} {\mathcal D}
\newcommand {\Z} {\mathcal Z}
\newcommand {\DD} {\mathbf D}
\newcommand {\bCP} {\mathbb {CP}^1}
\newcommand{\Pkt}{\textrm{.}}
\newcommand{\Kma}{\textrm{,}}
\newcommand {\Ga} {\Gamma}
\newcommand {\ga} {\gamma}
\newcommand {\eps} {\epsilon}
\newcommand {\de} {\delta}
\newcommand {\supp} {\mathrm supp~}
\newcommand {\HH} {\mathcal H}
\newcommand {\CC} {\mathcal C} 
\newcommand {\ST} {\mathcal {ST}} 
\newcommand {\SD} {\mathcal {SD}} 

\newcommand {\one} 
{\begin{picture}(12,8)
 \put(6,3.5){\circle{10}}
 \put(3.2,0.7){\text{1}}
\end{picture}
}
\newcommand {\oo}[1]
{\begin{picture}(12,8)
 \put(6,3.5){\circle{12}}
 \put(2.2,0.7){\text{$#1$}}
\end{picture}
}
\newcommand {\ooi}[1] 
{\begin{picture}(13,8)
 \put(6,3.5){\circle{13}}
 \put(1,0.7){\text{$#1$}}
\end{picture}
}
%%%%%%%%%%%%%%
%%% Another possibility $\makebox[0pt][l]{$\bigcirc$}\ 1$
%%%
%%%%%%%%%%%%%%%

\numberwithin{equation}{section}

             \title[Using homological duality in consecutive pattern avoidance]
             {Using homological duality in  consecutive pattern avoidance}
             
\author[A.~Khoroshkin]{Anton Khoroshkin}
\address{Departement Matematik, ETH, R\"amistrasse 101, 8092 Zurich, Switzerland and  ITEP, Bolshaya Cheremushkinskaya 25, 117259, Moscow, Russia} 
\email{anton.khoroshkin@math.ethz.ch}

\author[B.~Shapiro]{Boris Shapiro}
\address{Department of Mathematics, Stockholm University, SE-106 91, Stockholm,
            Sweden}
\email{shapiro@math.su.se}

\date{\today}

\begin{abstract}
Using the approach  suggested in \cite{DKh} we present below a sufficient condition guaranteeing  
that  two collections of patterns of permutations have the same exponential generating functions 
for the number of permutations avoiding elements of these collections  as consecutive patterns.  
In short,   the coincidence  of the latter generating functions is guaranteed by a  length-preserving 
bijection of  patterns in  these collections  which is identical on the overlappings of pairs of 
patterns where the overlappings are considered as unordered sets.
Our   proof is based on  a direct algorithm for the computation of the inverse generating functions.
As an application we present a large class of patterns where this algorithm is  fast and, in particular,  allows to obtain a  linear ordinary differential equation with polynomial coefficients 
satisfied by the  inverse generating function.
\end{abstract}

\maketitle

\section {Introduction} 

In the recent years the theory of consecutive pattern avoidance for permutations has experienced a rapid  development since the publication of the important paper \cite{ElN}.  Among the latest publications one should mention e.g. \cite{AAM}, \cite{LR}, \cite {Ra}, \cite{EKP} �where a  number of special cases was treated and the corresponding exponential generating function was explicitly found. The present  text is devoted to the same topic and is an extension of  the application of  homological methods to this theory  initiated in \cite{DKh}.  We investigate  a natural  analog of the notion of Wilf equivalence for consecutive pattern avoidance and obtain  a rather general sufficient condition  guaranteeing that this natural analog of Wilf equivalence holds. Most of the  definitions below  are borrowed from  \cite{DKh} and are (hopefully)  standard in this area. 

\subsection{Notation and definitions}

 A {\it permutation of length $n$} is a sequence containing each of the numbers $\{1,\ldots,n\}$ exactly once. To every sequence $s$  consisting of $n$ distinct positive integers, we associate  its {\it standardization $\st[s]$}  which is the permutation  of length $n$  uniquely determined by the condition that $s(i)< s(j)$ if and only if $\st[s](i)<\st[s](j)$. In other words, $\st[s]$ is the unique permutation of length $n$ whose relative order of entries is the same as that of~$s$.  For example, $\st[573]=(231).$ In what follows we will refer to  separate integers forming a permutation as its  {\em entries}.

We say that a permutation $\sigma$ of length $n$ {\it contains a permutation $\pi$ of length $l\le n$ as a consecutive pattern} if for some 
$i\leq n-l+1 $ the standardization $\st[\sigma(i)\sigma({i+1})\ldots\sigma({i+l-1})]$ coincide with $\pi$. If $\sigma$ contains $\pi$ as a consecutive pattern we say that 
{\it $\pi$ divides $\sigma$} and use the notation $\pi | \sigma$.  %we are going to use for the defined above notion of containing the subword will not mislead the reader since the same is frequently used in the theory of monomials in a free associative algebra. For example 
If  $\pi | \sigma$ and   $i=1$ (resp. $i=n-l+1$) we say that {\it $\pi$ is a left (resp. right) divisor of $\sigma$}. 
The main notion in the theory of pattern avoidance for permutations  is as follows. We say that {\it a  permutation $\sigma$ avoids a given permutation $\pi$ as a consecutive pattern} if $\sigma$ is not divisible by $\pi$.
 (Throughout this paper we only deal with consecutive patterns, so the word ``consecutive'' will be omitted.) 

The central problem of the theory of pattern avoidance is to count the number of  permutations of a given length avoiding a given collection~$\varPi$  of forbidden patterns or, more generally, containing a given number of occurrences of patterns from~$\varPi$. This problem naturally leads to the following equivalence relation on collections of patterns defined in the simplest case by H.~Wilf,~\cite{Wilf}. Two collections of patterns $\varPi_1$ and $\varPi_2$ are said to be {\it  Wilf equivalent} %(notation: $\varPi_1\simeq_W \varPi_2$) 
if for every positive integer  $n$, the number of $\varPi_1$-avoiding permutations of length~$n$ is equal to the number of $\varPi_2$-avoiding permutations of length~$n$. The next notion although very natural seems to be new. We say that two collections of patterns $\varPi_1$ and $\varPi_2$ are  {\it strongly Wilf equivalent} %(notation: $\varPi_1\simeq_{SW} \varPi_2$) 
if for every positive integer $n$ and every nonnegative integer $0\le q<n$, the number of permutations of length~$n$ with~$q$ occurrences of patterns from $\varPi_1$ equals  the number of permutations of length~$n$ with $q$ occurrences of patterns from~$\varPi_2$. However, 
standard Wilf equivalence deals with all pattern occurrences (and not just consecutive ones). In the set-up of consecutive pattern avoidance we will speak about {\it c-Wilf equivalent}� (resp.  {\it strongly c-Wilf equivalent})  collections where "c" stands for consecutive. (We  use the notation: $\varPi_1\simeq_{c-W} \varPi_2$ for strongly c-Wilf equivalent collections.)
\begin{remark}
 Throughout this paper we assume that   every collection of patterns $\varPi$ is {\it reduced}, i.e. no two permutations $\pi,\pi'\in \varPi$ are divisible by one another. (Notice that  if $\pi|\pi'\in \varPi$ then $\varPi\setminus\{\pi'\}$ is strongly c-Wilf equivalent to $\varPi$.) 
\end{remark}
Following \cite{ElN} consider    two natural exponential generating functions in one and two variables respectively:
$$\varPi(x):=\sum_{n}\al_{n}\frac{x^n}{n!}
 \quad \text{ and } \quad
\varPi(x,t):=\sum_{n,k}\al_{n,q}\frac{x^n}{n!}t^q,$$
associated to a given collection of patterns $\varPi$. 
 Here $\al_n$ (resp. $\al_{n,q}$)  is the number of permutations of length~$n$
avoiding all patterns from $\varPi$ (resp. the number of permutations of length $n$ with exactly ~$q$ occurrences of patterns from $\varPi$). Obviously, 
$\varPi(x)=\varPi(x,0)$. (Recall that we only count consecutive occurrences.) 

\begin{remark}� 
Hilbert series very similar to $\varPi(x)$ and $\varPi(x,t)$ are often considered in the theory of associative algebras.  The well-known method of their study is based on the so-called bar-cobar duality which roughly means that  a graded associative algebra  $A$ and the $A_{\infty}$-coalgebra $Tor^{A}_{\ldot}(\kk,\kk)$ are dual with respect to the functor $Tor$.  As a corollary of this duality  one gets that  the  Hilbert series of $A$ and of    $Tor^{A}_{\ldot}(\kk,\kk)$  are the inverses of  each other, i.e. their product equals $1$. (See \cite{Uf} for the details on different computational methods for the Hilbert series of associative algebras and their homology.)
It seems highly plausible that for an associative algebra with few relations  
a combinatorial description of its homology is simpler than that of the algebra itself. 
However, for algebras with many relations the situation is the opposite one. 
\end{remark} 

 Recall that the set of permutations avoiding an  arbitrary fixed collection $\varPi$ has an important additional structure, see appendix in \cite{DKh}. Namely, in  monoidal category  it can  be considered as the monomial basis of an algebra with monomial relations.  (We refer an interested reader to the above appendix in \cite{DKh} and references therein for the details. In particular, one can find  the definition of the homology functor in the latter appendix.)   
%One of the definitions of this monoidal structure may be found in \cite{Ro}. Corresponding algebras are called there {\em shuffle algebra}. 
Therefore, it seems natural to use the above mentioned homological duality in the theory of pattern avoidance.
 Combinatorial data appearing in this context is based on a generalization of  the so-called cluster method of I.~Goulden and D.~Jackson, \cite{GJ}. %Most of details may be found in \cite{DKh} as well.
We   explain below how one can get combinatorial information (for example, about the coefficients of the generating functions) of  the corresponding graded homological vectorspaces   for collections of patterns with few entries. 

%NADO OPREDELIT* OVERLAPPING!!!

To describe our results we need to  recall the  definition of a combinatorial gadget called {\it clusters} in \cite{GJ}. 
They generalize the notion of a linkage given below.  

A permutation $\sigma$ of length $n$ is called a {\it linkage} of an ordered pair of (not necessarily distinct) patterns $(\pi, \pi')$ of lengths $l$ and $l'$ if (i)  $n<l+l'$; and (ii) the standardizations $\st[\sigma(1)\ldots\sigma(l)]$
and $\st[\sigma(n-l'+1)\ldots\sigma(n)]$ are equal to $\pi$ and $\pi'$ resp. 
Since the length of  $\sigma$ is less than the sum of the lengths of $\pi$ and $\pi'$ one has that the standardizations of the right truncation of $\pi$ and left truncations of $\pi'$ should be the same.
Setting  $k=(l+l'-n)$  we say that a pair $(\pi,\pi')$  has a {\it $k$-overlapping} or {\it $k$-overlaps}.
In other words,  a pair $(\pi, \pi')$ $k$-overlaps whenever the standardization $st[\pi(l-k+1)\ldots\pi(l)]$ is equal to the standardization $st[\pi'(1)\ldots\pi'(k)]$.
(Notice that  there might be several different linkages of two given patterns $\pi$ and $\pi'$ of the same length.)  

%%%%%%%%%%%%
%%%%%%%%%%%%%%%
%It is also relevant to have in mind the following picture:
%%%%%%%%%%%
%%%%Нарисовать картинку с зацеплением и определить что такое overlapping
%%%%%%%%%%

 A cluster is a way to link together several patterns from a given set. 
More precisely, % in the case of permutations
a $q$-cluster w.r.t. a  given collection of patterns $\varPi$ is a triple $(\sigma;\pi_1,\ldots,\pi_q;d_1,\ldots,d_q)$ where $\sigma$ is a permutation, $\{\pi_i\}$ is a list of (not necessarily distinct) patterns from $\varPi$, and $\{d_i\}$ is a list of positive integers such that
\begin{itemize}
\item[(i)] for every $j=1,\ldots,q$, $\st[\sigma({d_j}),\ldots,\sigma({d_j+l_j-1})]=\pi_j\in \varPi$, where $l_j$ is the length of~$\pi_j$ (here $d_j$ labels  the beginning  of the  pattern $\pi_j$ in $\sigma$);
\item[(ii)] $d_{j+1}>d_j$ (patterns are listed from left to right) and $d_{j+1}<d_j+l_j$ (adjacent patterns are linked);
\item[(iii)] $d_1=1$, and the length of $\sigma$ is equal to $d_q+l_q-1$ (i.e. $\sigma$ is completely covered by the patterns $\pi_1,\ldots,\pi_q$).
\end{itemize}
\medskip
Denote by $cl_{n,q}(\varPi)$  the number of $q$-clusters of length~$n$ in a collection $\varPi$ and introduce the exponential generating function
 $$\varPi_{cl}(x,t)=x+\sum_{n>1,q\geq 1}cl_{n,q}\frac{x^n}{n!}t^q.$$ 
(Here we use a natural convention  that there always exists  exactly one (fictitious) $0$-cluster and, therefore,  the above generating function starts with $x$.) 

The next  result is an immediate consequence of the general cluster method of I.~Goulden and D.~Jackson, \cite{GJ} and its homological proof for the case of permutations can  be found in \cite{DKh}. 

\begin{theorem}\label{th:GJ}
In the above notation one has: 
\begin{equation}\label{enum_cluster}
\varPi(x,t)=\frac{1}{1-\varPi_{cl}(x,t-1)}. 
\end{equation}
\end{theorem}

\begin{corollary} 
The exponential generating function $\varPi(x):=\varPi(x,0)$ of the numbers of permutations avoiding  patterns from a given  collection $\varPi$ satisfies the relation:  
\begin{equation}\label{avoid_cluster}
\varPi(x)=\frac{1}{1-\varPi_{cl}(x,-1)}.
\end{equation}
\end{corollary} 

\medskip 
\begin{remark} 
In general, the problem of counting the number of $q$-clusters in a given collection of patterns $\varPi$ does not seem to be easier than counting the number of permutations of a given length avoiding $\varPi$.  On the other hand, there exist natural classes of collections   for which counting $q$-clusters is an easier task, see \S~3.
\end{remark} 

One can guess that since clusters can  be described in terms of linkages of pairs of patterns the number of clusters can also be determined in terms of  the combinatorics of these linkages. Exploiting the latter idea we were able to prove the following (main) result of this note.  

\begin{theorem}\label{th:main}
Two collection of patterns $\varPi_1$ and $\varPi_2$ are strongly c-Wilf equivalent if there exists a bijection $\varphi: \varPi_1\rightarrow \varPi_2$ preserving  the following  three properties:
\begin{itemize}
 \item {\rm (lengths):} For any $\pi\in\varPi_1$ its length equals to that of $\varphi(\pi)\in \varPi_2$; 
 \item {\rm (linkages):} A pair of  patterns $(\pi, \pi')$ from $\varPi_1$ has a linkage of length $n$ if and only if the pair of its images  $(\varphi(\pi),\varphi(\pi'))$ from $\varPi_2$ has a linkage of the same length~$n$.
 \item {\rm (overlapping sets):} For each overlapping of any pair of  patterns from  $\varPi_1$ 
  the bijection $\varphi$ preserves  the subsets of entries that overlaps. More precisely,  for any pair  $(\pi,\pi')$ of patterns $\pi,\pi'\in\varPi_1$ of lengths $l$ and $l'$ resp.  and an arbitrary positive integer  $k\leq min(l,l')$ the 
coincidence of the standardizations 
$\st[(\pi(l-k+1)\ldots\pi(l))] = \st[(\pi'(1)\ldots\pi'(k))]$
 implies the coincidence of sets:  
$$ \{\pi(l-k+1),\ldots,\pi(l)\}=\{\varphi(\pi)(l-k+1),\ldots,\varphi(\pi)(l)\}, \text{   and   } 
\{\pi'(1),\ldots,\pi'(k)\}=\{\varphi(\pi')(1)\ldots\varphi(\pi')(k)\}.$$
\end{itemize}
\label{theorem_Wilf_overaping}
\end{theorem}

%NUZHEN PRIMER 2 NABOROV POZAMYSLOVATEE. 

The simplest case where Theorem~\ref{th:main} applies is to  collections with a single pattern having no nontrivial self-overlappings.
The following result implied by Theorem~\ref{th:main}   was first conjectured by S.~Elisalde in \cite{El} and  later proved  in \cite{DKh} by homological methods and, simultaneously,  by J.~Remmel whose methods were based on \cite{MR}. 
Namely, 
\begin{corollary}
 Two collections of patterns each containing a single permutation  without nontrivial self-overlappings are strongly c-Wilf equivalent if 
 \begin{itemize}
 \item[(i)] the lengths of the permutations coincide;
  \item[(ii)] the first entry and resp. the last entry of the permutations coincide.
 \end{itemize}
\end{corollary}

A series of particular examples covered by Theorem~\ref{th:main} can  be found in \S~5 of \cite{AAM}. 
Namely, the following definitions are borrowed from the latter paper.

We say that a pair of permutations $\alpha\in S_k$ and $\beta\in S_{k'}$ {\it has a separation property} if
$\beta$ avoids the permutation $(\alpha(1),\ldots,\alpha(k),k+1)\in S_{k+1}$ 
and $\alpha$ avoids $(k'+1,\beta(1),\ldots,\beta(k'))\in S_{k'+1}$.

With each pair of permutations $\alpha\in S_k$, $\beta\in S_k'$ and a natural number $l$
one can associate the subset $\Pi(\alpha,\beta;l)\subset S_{k+l+k'}$ of permutations  defined by two following properties.
We say that $\pi\in\Pi(\alpha,\beta;l)$ iff
\begin{itemize}
 \item[(i)] the standardizations of the $k$ first and $k'$ last entries coincide with $\alpha$, and $\beta$ resp.
 \item[(ii)] the $k$ first entries are strictly smaller than the $k'$ last entries; 
the $k'$ last entries are strictly smaller than the remaining entries of $\pi$ in the middle.
In other words, $\pi(i)<\pi(j)<\pi(s)$ for any triple of indices $(i,j,s)$ such that $1\leq i\leq k < s \leq k+l< j\leq k+l+k'$.
\end{itemize}
\begin{corollary}
Fix a pair of permutations  $\alpha$ and $\beta$ having a separation property and a
 $d$-tuple of natural numbers $(l_1,\ldots,l_d)$. 
Then all collections of $d$ distinct patterns $\{\pi_1,\ldots,\pi_d\}$ such that $\pi_i\in\Pi(\alpha,\beta;l_i)$ are strongly Wilf equivalent.
\end{corollary}
\begin{proof}
 For all $i$ the $l_i$ entries in the middle of each pattern $\pi_i\in\Pi(\alpha,\beta,l_i)$ never appear in the overlapping sets.
\end{proof}
%%%%%%%%%%%%%%

Let us present a few  more examples illustrating how our theorem works in practice.

The following  patterns 
 $$
(1734526)\sim_{\text{c-W}} (1735426) \sim_{\text{c-W}} (1743526) \sim_{\text{c-W}}
(1745326) \sim_{\text{c-W}} (1753426) \sim_{\text{c-W}} (1754326) 
$$
are pairwise c-Wilf equivalent. 
They have self-overlappings of lengths $1$ and $2$ and coinciding pairs of the 
first two and the last two entries. 
 
The next pair of c-Wilf equivalent patterns 
\begin{equation}
\label{wilf_example}
(143265987) \sim_{\text{c-W}} (134265897)
\end{equation}
have self-overlappings of lengths $1$ and $4$ and the corresponding 
subsets of their initial and final entries of lengths 1 and 4 coinciding while  their initial and final subwords are different. 

Finally, here is an example  
$$\{(145623),(13452)\} \sim_{c-W} 
\{(145623),(13542)\} \sim_{c-W}
\{(146523),(13452)\} \sim_{c-W}
\{(146523),(13542)\}
$$
of c-Wilf equivalent collections with 2 patterns in each.

In \S~2 we prove Theorem~\ref{th:main} and in \S~3 we apply our main construction to a class of collections of patterns and obtain a system of linear ordinary differential equations satisfied by $\Pi_{cl}(x,t)$ together with a set of similar generating functions defined below.  In the follow-up \cite{KhS}  of the present  paper we plan to study different asymptotic properties of $\Pi(x,t)$ using the suggested approach. 

\medskip 
\noindent
{\it Acknowledgements.}  The authors are sincerely grateful to  S.~Kitaev for the  e-mail correspondence  concerning this subject. 
Research of the first author is supported by RFBR 10-01-00836, RFBR-CNRS-10-01-93111, RFBR-CNRS-10-01-93113,
and by the Federal Programm under the contract 14.740.11.0081.

\section {Proofs} 
\label{sec:proofs}

Our proof of Theorem~\ref{theorem_Wilf_overaping} consists of  an algorithm  computing the cluster generating function $\varPi_{cl}(x,t)$ of a given collection of patterns $\varPi$. It will  then be relatively easy  to see that  this algorithm uses only the lengths and the overlapping subwords for pairs of  patterns from $\varPi$ considered as sets. To start with,  we define 
for an arbitrary  collection of patterns  $\varPi$  a certain directed graph with labelled vertices and edges. 
The important  circumstance  is that the number of $q$-clusters with fixed initial and final subwords  will be equal to the number of properly 
weighted paths of length $q$ in this graph with fixed initial and final vertices. The required weights can be computed using the edge labels. 
As a consequence  we get that  this graph  uniquely defines the generating functions  $\varPi_{cl}(x,t)$ and, therefore,  $\varPi(x,t)$, see Theorem~\ref{th:GJ}. 

Namely, given an arbitrary collection of patterns $\varPi$ define its directed  graph $\mathcal{G}(\varPi)$ with labelled vertices and edges as follows. The vertices of 
$\mathcal{G}(\varPi)$ will be labelled by permutations (of, in general, different lengths) and the labels of the edges are defined below. 
\begin{itemize}
 \item To define the vertices assume that some permutation ${v}$ is a left divisor  of a pattern $\pi_{\alpha}\in \varPi$ and, 
at the same time,  a right divisor of a (not necessarily different) pattern $\pi_{\beta}\in \varPi$. 
Then we assign to ${v}$ a vertex $\oo{v}$ of $\mathcal{G}(\varPi)$ and, naturally,  label this vertex by $v$. 
(Notice that the same $v$ can arise from different pairs $(\pi_{\alpha},\pi_{\beta})$. 
In particular,   the trivial  $1$-element permutation $(1)$ comes from an arbitrary pair of not necessarily distinct patterns.  
$\one$ is called the {\it distinguished vertex} of $\mathcal{G}(\varPi)$ and  the set of all vertices of  
$\mathcal{G}(\varPi)$ is denoted by $\mathcal{V}(\varPi)\ni \one$.)
 
\item To define the edges take a pattern $\pi\in\varPi$ of some length $l$ and a pair $(\pi_i,\pi_f)$ of its initial and final subwords of lengths $k$ and $k'$ (i.e. $\pi_{i}:=(\pi(1)\ldots\pi(k))$ and $\pi_f:=(\pi(l-k'+1)\ldots\pi(l))$) such that 
standardizations $\st[\pi_i],\st[\pi_f]$ are the vertices of  $\mathcal{G}(\varPi)$.  %, i.e. belong to $\mathcal{V}(\varPi)$. 
 Let $\mu_i$ and $\mu_f$ be the subsets of entries that appears in $\pi_i$ and $\pi_f$ resp. (i.e. $\mu_{i}:=\{\pi(1),\ldots,\pi(k)\}$ and $\mu_f:=\{\pi(l-k'+1),\ldots,\pi(l)\}$).
 The triple $(\pi, \pi_i, \pi_f)$ then defines the directed edge from the vertex $\st[\pi_i]$ to the vertex $\st[\pi_f]$ which we label   by the triple $(\mu_i,\mu_f;l)$.
\end{itemize}

\begin{remark} Notice that $\mu_i$ and $\mu_f$ are considered as {\it unordered sets}. 
\end{remark} 

\noindent
{\it Notation.} In the present  text  we show all vertices of $\mathcal{G}(\varPi)$ as encircled permutations to distinguish them from their labeling  permutations. 
We will try to keep our notation straight throughout the whole text by denoting  similar quantities by the same letter and   adding extra indices if required. 
For example,  $l$ will typically mean the length of a pattern $\pi$ from a collection, 
$k$ should denote the length of a permutation $v$ which labels a vertex  of $\mathcal{G}(\varPi)$ originating  from a $k$-overlapping, $n$ will stand for the  length of a cluster.

Four examples of  $\mathcal{G}(\varPi)$ are given below.  
The upper left example is constructed from the  collection $\varPi_1 =\{(1342765),(152364)\}$ of two patterns with no nontrivial overlappings. The upper right example comes from the single pattern 
$\{(132679485)\}$ having self-overlappings of lengths $1$ and $3$. The meaning of two other examples should be clear now. 
\begin{figure}[htb]
\includegraphics{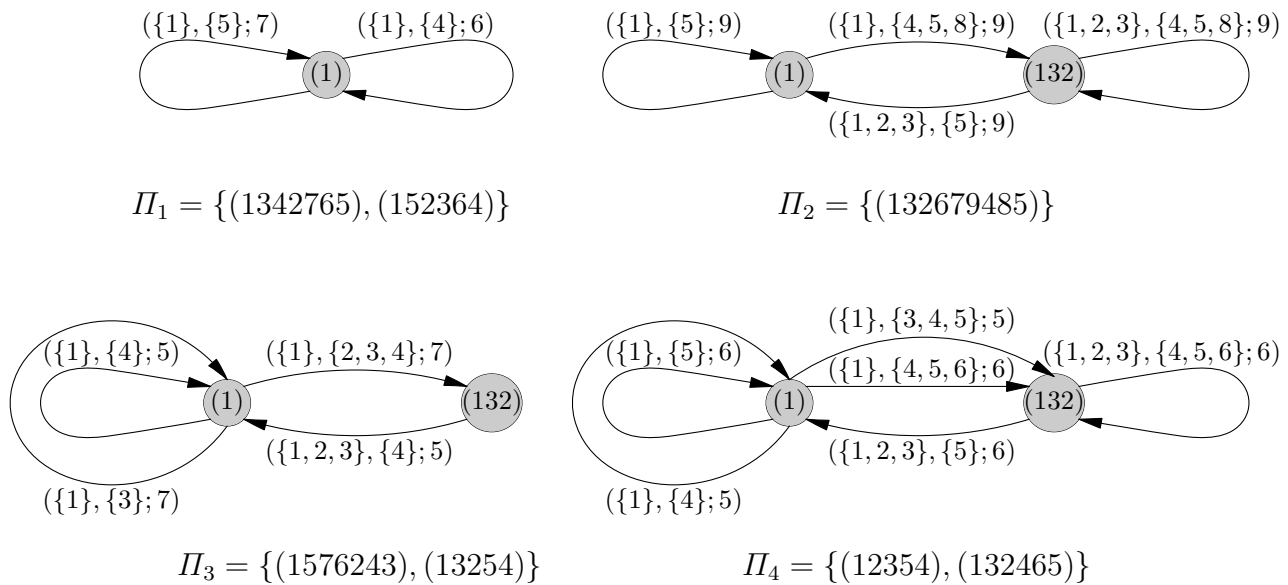}
\caption{Four examples of $\mathcal{G}(\varPi)$.}
\label{graph-example}
\end{figure}

  Our main technical result  is as follows. 

\begin{theorem}\label{th:graph} 
The graph $\mathcal{G}(\varPi)$  uniquely determines the generating function $\varPi_{cl}(x,t)$. 
\end{theorem}
The next corollary immediately implies Theorem~\ref{theorem_Wilf_overaping}.
\begin{corollary} Two  collections of patterns $\varPi_1$ and $\varPi_2$ having isomorphic graphs $\mathcal{G}(\varPi_1)$  and $\mathcal{G}(\varPi_2)$ are strongly c-Wilf-equivalent. (Here by an 'isomorphism'  we mean a  graph isomorphism preserving the labels of edges. The labels of vertices can change.) 
\end{corollary}

\begin{proof}
\label{algorithm}
To prove Theorem~\ref{th:graph} we  present a natural algorithm calculating the number of 
$q$-clusters in a given collection $\varPi$ in terms of its  graph $\mathcal{G}(\varPi)$. 
Namely, each vertex $\oo{v}$ and a positive integer  $n$ uniquely determine the subset $Cl_{v,n,q}$ consisting of all 
$q$-clusters $(\sigma;\pi_1,\ldots,\pi_q;d_1,\ldots,d_q)$, such that 
the length of $\sigma$ is equal to $n$ and the standardization of the initial subword of $\sigma$ is equal to $v$. 
Moreover,   with each word $\bar{p}:=(p_1\ldots p_{k})$ of length $k$ (where $k$ is the length of  $v$) 
one can associate the subset $Cl_{v,n,q}[\bar{p}]\subset Cl_{v,n,q}$ consisting of those clusters in $Cl_{v,n,q}$ 
which  have $\bar{p}$ as their initial subword.
We will explain how one should compute the cardinalities of  $Cl_{v,n,q}[\bar{p}]$ by induction on $q$
using the edge labels in  $\mathcal{G}(\varPi)$. 
Therefore, the cardinalities of $Cl_{v,n,q}$  can also be computed inductively as the sums over different $\bar{p}$.
Since the standardization of any word of length $1$  equals $(1)$  the set $Cl_{(1),n,q}$ 
coincides with the set of all $q$-clusters of length $n$. (The cardinality of the latter set is one of required coefficients in the cluster generating function $\varPi_{cl}(x,t)$.)

Let us now return to  the induction step. Take an arbitrary vertex $\oo{v}\in \mathcal{V}(\varPi)$ and  
let $\oo{v}\stackrel{\pi_1}{\mapsto} \ooi{v_1},$\ldots, 
$\oo{v}\stackrel{\pi_d}{\mapsto} \ooi{v_d}$ 
 be the list of all edges in  $\mathcal{G}(\varPi)$ starting at the vertex $\oo{v}$.  
Denote by $k_j$  the length of permutation $v_j$ labeling the  vertex $\ooi{v_j}$  and 
denote by  $l_j$ the length of the pattern $\pi_j$. 
We present below a recurrence  relation expressing the cardinality $cl_{v,n,q}[\bar{p}]$  of the set $Cl_{v,n,q}[\bar{p}]$ 
in terms of the cardinalities  $cl_{v_j,n-l_j+k_j,q-1}[\bar{p}']$ of $Cl_{v_j,n-l_j+k_j,q-1}[\bar{p}']$
with the summation  taken over a certain subset of words $\bar{p}'$. 
Using this relation we can calculate inductively  each $cl_{v,n,q}[\bar{p}]$ and then obtain 
the required $cl_{v,n,q}$ by summation over different $\bar{p}$.
It will be  convenient to subdivide the sets $Cl_{v,n,q}$ and $Cl_{v,n,q}[\bar{p}]$ into the subsets indexed by the  edges starting at  the vertex $\oo{v}$. For example, $Cl_{v\stackrel{\pi_j}{\mapsto} v_j,n,q}$ is the subset of $q$-clusters formed by  linkages of length $n$ between  the pattern $\pi_j$ and a $(q-1)$-cluster from $Cl_{v_j,n-l_j+k_j,q-1}$. One has 
\begin{equation}
cl_{v,n,q} = 
\sum_{
\begin{smallmatrix}
 { 1\leq p_1,\ldots,p_k \leq n, }\\
  st[(p_1\ldots p_k)]=v
\end{smallmatrix}
} cl_{v,n,q}[p_1\ldots p_k] =
\sum_{
\begin{smallmatrix}
 { 1\leq p_1,\ldots,p_k \leq n, }\\
  st[(p_1\ldots p_k)]=v
\end{smallmatrix}
} 
 \sum_{j=1}^d cl_{v\stackrel{\pi_j}{\mapsto} v_j,n,q}[p_1\ldots p_k]. 
\label{add_statistics}
\end{equation}

Therefore,  it is sufficient  to find  recurrence relations   expressing the terms $cl_{v\stackrel{\pi_j}{\mapsto} v_j,n,q}[p_1\ldots p_k]$   
in  the right-hand side of \eqref{add_statistics} using $cl_{v_j,n-l_j+k_j,q-1}[...]$.
To avoid very cumbersome  notation let us take the case of a single edge starting at $\oo{v}$ which is equivalent to fixing $v_j$ in the above formulas. 
Let  $\oo{v}\stackrel{\pi}{\mapsto}\ooi{v'}$ be an edge in a graph $\mathcal{G}(\varPi)$
coming from a pattern $\pi$ of length $l$ and 
let $k$ and $k'$ be the lengths of the permutations labeling $\oo{v}$ and $\oo{v'}$ resp. 
To explain our  recurrence we need to introduce the following extra notation associated to $\pi$. 

If $l>k+k'$ let $\psi\in S_{k+k'}$ be the permutation which is the inverse of the standardization 
of the $k$ first  and the $k'$ last  entries of $\pi$, and let $\overline{\psi}$ be the 
composition of $\psi$ with the  shifting map 
$sh_{k,k'\mapsto l}:\{1,\ldots,k,k+1,\ldots,k+k'\}\rightarrow \{1,\ldots,k\}\cup\{l-k'+1,\ldots,l\}$ 
defined by the formula:
$$ sh_{k,k'\mapsto l}(j) = 
\left\{
\begin{array}{l}
 j,\text{ if } j\leq k, \\
 j+ l-k-k'+1, \text{ if } j > k.
\end{array}
\right.
$$
In other words, $\overline{\psi}$ prescribes the rule how to write down the $k$ first  and the $k'$ last  
entries of the pattern $\pi$ in the increasing order:
$$
 \{\pi(\overline{\psi}(1)) <\pi(\overline{\psi}(2)) <\ldots < \pi(\overline{\psi}(k+k'))\} = 
\{\pi(1),\ldots,\pi(k)\}\cup\{\pi(l-k'+1),\ldots,\pi(l)\}. 
$$

The next statement gives the required recurrence. 

\begin{lemma}\label{lm:1}
The following relations hold: 
\begin{itemize}
\item
 for $l>k+k'$ set $\tilde{\pi}=st[\pi(1)\ldots\pi(k)\pi(l-k'+1)\ldots\pi(l)]$. Then 
\begin{multline} 
cl_{v\stackrel{\pi}{\mapsto} v',n,q}[p_{1}\ldots p_{k}]=
 \\
\sum_{
\begin{smallmatrix}
 p_{k+1},\ldots,p_{k+k'}\ :\\
st[(p_1\ldots p_{k+k'})]= \tilde{\pi}
\end{smallmatrix}}
\binom{p_{\psi(1)}-1}{\pi(\overline{\psi}(1))-1}\times 
\left[
\prod_{j=1}^{k+k'-1}
\binom{p_{\psi(j+1)}-p_{\psi(j)}-1}{\pi(\overline{\psi}(j+1))-\pi(\overline{\psi}(j))-1}
\right]
\times \binom{n-p_{\psi(k+k')}}{l-\pi(\overline{\psi}(k+k'))}\times \\
\times
cl_{v',n-l+k',q-1}[p_{k+1}-\pi(l-k'+1) +v'(1),\ldots,p_{k+k'}-\pi(l)+v'(k')] .
\label{one_edge}
\end{multline}
\item 
for $l\leq k+k'$ one has: 
\begin{multline}
{ cl_{v\stackrel{\pi}{\mapsto}v',n,q}[p_{1}\ldots p_{k}]=
} 
\\ 
{ =
\sum_{p_{k+1},\ldots,p_{l} : \st[(p_1\ldots p_l)]= \pi} 
cl_{v',n-l+k',q-1}[p_{l-k'+1}-\pi(l-k'+1) +v'(1),\ldots,p_{l}-\pi(l)+v'(k')] .}
\label{one_edge_1}
\end{multline}
\end{itemize}
\end{lemma}

\begin{remark} 
The range of  summation in \eqref{one_edge} can be easily derived from  our  convention 
on the binomial coefficients claiming that $\binom{N}{M}=0$ if either $N<0$ or $M>N$. 
Moreover, we assume that  $p_j$'s are pairwise different positive integers not  exceeding $n$.
For  the induction base we use the following initial data:
$$
Cl_{v,n,0} = \left\{
\begin{array}{l}
 {\{(1)\} \text{, if } v=(1) \text{ and } n=1, }\\ 
{ \varnothing,  \text{ otherwise. } }
\end{array}
\right.
$$
\end{remark} 

\begin{proof}
We show how  to prove  \eqref{one_edge}. 
In formula \eqref{one_edge} one has the summation over all patterns $\sigma\in Cl_{v\stackrel{\pi}{\mapsto} v',n,q}$ such that the 
word $(\sigma(1)\ldots\sigma(k)\sigma(l-k'+1)\ldots\sigma(l))$ is fixed and coincides with $(p_1\ldots p_{k+k'})$.
Indeed, the numbers $p_j$ are ordered by the permutation $\psi$ as follows: $p_{\psi(1)}<\ldots<p_{\psi(k+k')}$. 
Therefore, there are $\binom{p_{\psi(1)}-1}{\pi(\overline{\psi}(1))-1}$ choices of entries less than $p_{\psi(1)}$ among the first  
$l$ entries of $\sigma$; there are 
$\binom{p_{\psi(2)}-p_{\psi(1)}-1}{\pi(\overline{\psi}(2))-\pi(\overline{\psi}(1))-1}$ choices of  entries greater than 
$p_{\psi(1)}$ 
and less than $p_{\psi(2)}$, $\ldots$;  there are  
$\binom{n-p_{\psi(k+k')}}{l-\pi(\overline{\psi}({k+k'}))}$ choices of  entries greater than $p_{\psi(k+k')}$
among the first $l$ entries of $\sigma$; and 
$cl_{v',n-l+k',q-1}[p_{k+1}-\pi(l-k'+1) + v'(1),\ldots,p_{k+k'}-\pi(l)+ v'(k')]$
ways to choose the remaining standardization of the last $(n-l+k')$ entries of $\sigma$.

In  \eqref{one_edge_1}   the union of the $k$ initial  entries and the $k'$ final  entries of $\pi$
covers  the whole list of entries of $\pi$, i.e. the set  $\{1,\ldots,l\}$. Therefore,  all binomial coefficients appearing  
in  \eqref{one_edge} should be equal to $1$ which leads to \eqref{one_edge_1}. 
\end{proof}
As an immediate consequence of Lemma~\ref{lm:1}  one can see that the numbers  
$cl_{v\stackrel{\pi}{\mapsto}v',n,q}[\ldots]$ of  $(q+1)$-clusters  depend only on the length, 
the $k$  first and the $k'$ last  entries of  $\pi$ considered as sets. 
This justifies  the information we use as the edge labels  of the graph $\mathcal{G}(\varPi)$.
The formulas expressing $cl_{v,n,q}[\ldots]$ in terms of $cl_{\ldot,\ldot,q-1}[\ldots]$  
depend only on the labeling  of the edges starting at $\oo{v}$.
Therefore,  these cardinalities can  be computed by induction on $q$
using the edge labels of the  graph $\mathcal{G}(\varPi)$.
Finally, as we mentioned before,  the set of all $q$-clusters of length $n$ of the whole collection $\varPi$ is equal to the set 
$Cl_{(1),n,q}$\,.  
\end{proof}

\subsection{Case of a single pattern}
\label{section-one-pattern}
Let us consider separately the situation when  $\varPi$ contains  just a single  pattern since  in this case some substantial simplifications of our construction can be done.

First of all the following observation explains why the graph $\mathcal{G}(\{\pi\})$ is not required.
\begin{lemma}
 Let $\pi$ be a pattern of length $l$ and let $(2l-k_1),$\ldots,$(2l-k_d)$ be the list of all distinct lengths 
of possible self-linkages of $\pi$, i.e. $k_1$, \ldots, $k_d$ is the list of distinct lengths of self-overlappings of $\pi$.
Then  $\mathcal{G}(\{\pi\})$ is a complete directed graph on  $d$ vertices with loops and with lengths of the underlying permutations being equal to $k_1,\ldots,k_d$.  
Each ordered pair of  (not necessary distinct) vertices  of $\mathcal{G}(\{\pi\})$ are connected by exactly one directed edge labeled by the corresponding initial and final subwords of $\pi$.
\label{one_pattern}
\end{lemma}
It is obvious that $k_1=1$ and  denote  by $k$ ($k=k_d$) the length of the largest  overlapping.
Let $v_s$ be the standardization of the $k_s$ first  entries of $\pi$ 
(i.e. $v_s$ is the labeling permutation of the $s$-th vertex in $\mathcal{G}(\{\pi\})$).
Since all patterns  involved in any cluster coincide with $\pi$ 
the standardization of the initial subword of any cluster should  always be the same. 
Hence  for different $v_s$ and fixed $n$ and $q$ all the sets $Cl_{v_s,n,q}$ coincide.
Therefore,  it makes sense  to denote by $cl_{n,q}$ and $cl_{n,q}[p_1\ldots p_k]$ the cardinalities of 
the set of $q$-clusters of length $n$ and those  having $(p_1\ldots p_k)$ as their initial subword resp.
We  introduce the same set of notation for the self-overlappings of $\pi$  similar to what  have been using in   Lemma~\ref{lm:1} for the case $k_s<l-k$.

Namely, for $l>k+k_s$ let $\psi\in S_{k+k_s}$ be the permutation which is the inverse of the standardization 
of the $k$ first  and the $k_s$ last  entries of $\pi$; for $l\leq k+k_s$ let $\psi$ be the inverse of $\pi$.
Let $\overline{\psi_s}$ be the composition $sh_{k,k_s\mapsto l}\circ \psi_s$ using which  one gets the following rearrangement 
of the first $k$ and last $k_s$ elements of $\pi$ in the increasing order: 
$$
 \{\pi(\overline{\psi_s}(1)) <\pi(\overline{\psi_s}(2)) <\ldots < \pi(\overline{\psi_s}(k+k_s))\} = 
\{\pi(1),\ldots,\pi(k)\}\cup\{\pi(l-k_s+1),\ldots,\pi(l)\}. 
$$
Additionally,  let $\tilde{\pi}_s$ be the standardization 
of the $k$ first  and $k_s$ last  entries of $\pi$.

 In the case of a single pattern Lemma~\ref{lm:1} implies the next result. 

\begin{lemma}\label{lm:2}
For a single pattern the recurrence formula for the numbers of $q$-clusters is as follows:
\begin{multline}
cl_{n,q}[p_1\ldots p_k] = \\
\sum_{
\begin{smallmatrix} 
       s \ :\\ 
      k_s < l-k       
\end{smallmatrix} 
}
\sum_{
\begin{smallmatrix}
 1\leq p_{k+1},\ldots,p_{2k}\leq n, \\
st[(p_1\ldots p_{k+k_s})]= \tilde{\pi}_s
\end{smallmatrix}}
\binom{p_{\psi_s(1)}-1}{\pi(\overline{\psi_s}(1))-1}\times
\left[ \prod_{j=1}^{k+k_s-1}   
\binom{p_{\psi_s(j+1)}-p_{\psi_s(j)}-1}{\pi(\overline{\psi_s}(j+1))-\pi(\overline{\psi_s}(j))-1} \right]  
\times \binom{n-p_{\psi_s(k+k_s)}}{l-\pi(\overline{\psi_s}({k+k_s}))}\times 
\\
\times 
cl_{n-l+k_s,q-1}[p_{k+1}-\pi(l-k_s+1) +v_s(1),\ldots,
p_s({k+k_s})-\pi(l)+v_s(k_s),p_{k+k_s+1},\ldots,p_{2k}]  
+ \\
+
\sum_{
\begin{smallmatrix} 
       s \ :\\ 
      k_s \geq l-k       
\end{smallmatrix} 
}  
\sum_{
\begin{smallmatrix}
 1\leq p_{k+1},\ldots,p_{l+k-k_s}\leq n, \\
st[(p_1\ldots p_{l})]= \pi
\end{smallmatrix}
}
cl_{n-l+k_s,q-1}[p_{l-k_s+1}-\pi(l-k_s+1) +v_s(1),\ldots,p_{l}-\pi(l)+v_s(k_s),p_{l+1},\ldots,p_{l+k-k_s}].
\label{recurent-one-pattern} 
\end{multline}
  
(As above we assume that the set of  $0$-clusters contains the unique fictitious element of length $1$ while the set of $1$-clusters contains the single pattern $\pi$.) 
\end{lemma}

\section{Application} 
In this section we  discuss a specific class  of collection of patterns. Our method from \S~\ref{sec:proofs}  allows us to  construct a system of linear ordinary differential equations in  the variable $x$ for the cluster  generating functions $\varPi_{cl}(x,t)$ together with a set of  similar generating functions associated to a given collection of patterns $\varPi$. The main definition is as follows. 

\begin{definition} A collection  of patterns $\varPi$ is called {\it monotone} if for all $k>0$ and for each pair (not necessarily distinct) patterns $(\pi,\pi')$ from  $\varPi$ the existence of their $k$-overlapping  implies that  the initial subword of the  pattern $\pi'$ does not contain entries greater than $k$. 
\end{definition}

The next lemma explains how the monotonicity assumption  simplifies the structure of the  set of clusters and their generating functions. %We hope that the word degenerate used in the title of the subject will also be clear from this point of view.
\begin{lemma}
\label{lm:monot_link}
 Let $\sigma$ be a linkage of a pair of  patterns $(\pi,\pi')$. 
Suppose that the initial subword of $\pi$ of length $k$ does not contain entries greater than $k$
and that the  initial subword of $\pi'$ of length $k'$ does not contain entries greater than $k'$ (where $k'$ is the length of the overlapping  of the pair  $(\pi,\pi')$ in $\sigma$).
Then the initial subword of length $k$ of $\sigma$ should be equal to the initial  subword of $\pi$,  i.e. $\sigma(j)=\pi(j)$ for $1\leq j\leq k$.
\end{lemma}
\begin{proof}
Let $l$ and $l'$ be the lengths of   $\pi$ and $\pi'$ resp. Then the  length of $\sigma$ should be equal to $n=l+l'-k'$. 
Set $j=\pi^{-1}(1)$ and $j'=\pi'^{-1}(1)$  respectively.
The  number $j$ is the minimal entry of $\pi$ and should be less than or equal to $k$
($j'\leq k'$ respectively).  Therefore, $\sigma(j)$  should be the smallest entry among the $l$ first  entries of the linkage $\sigma$ and $\sigma(n-l'+j')$ should be the smallest entry among  the $l'$ last  entries of $\sigma$. In particular,  $\sigma(j)<\sigma(n-l'+j')$.  This implies that $\sigma(j)$ is the smallest entry in the whole $\sigma$ and hence $\sigma(j)=\pi(j)=1$. 
Similar arguments based on simple properties of   standardizations  imply that
$\sigma(\pi^{-1}(2))$ is the only entry in $\sigma$ greater than $1$  etc.
\end{proof}

Lemma \ref{lm:monot_link} implies that for an arbitrary vertex $\oo{v}\in\mathcal{V}(\varPi)$ 
 all clusters in the subset  $Cl_{v,n,q}$ should always have $v$ as their initial subword.
This means that the number $cl_{v,n,q}[p_1\ldots p_k]$ is non-vanishing  if and only if 
$(p_1\ldots p_k)=(v(1)\ldots v(k))$, where $k$ is the length of $v$. 
Consider two vertices $\oo{v}$ and $\ooi{v'}$ of lengths $k$ and $k'$  connected by an edge 
$\oo{v} \stackrel{\pi}{\mapsto} \ooi{v'}$ in $\mathcal{G}(\varPi)$.
We can simplify  the formula \eqref{one_edge} using  our assumption on the 
intersections of clusters. At first consider  the case when the length $l$ of the pattern $\pi$ is greater than $k+k'$.
One observes that in this case  there is no summation in   \eqref{one_edge}  since the summands in the 
left-hand  and the right-hand  sides of  \eqref{one_edge}  are non-vanishing if and only if for $1\leq j\leq k$ 
one has that $p_{j}= v(j) = \pi(j)$ and for  $1\leq j\leq k'$ one has that $p_{k+j} = \pi(l-k'+j).$  
Therefore, only the last binomial coefficient in the product in the right-hand side of \eqref{one_edge}  is different from $1$.  
Denote  by $m$ the maximal entry of the final subword of $\pi$. 
(In notation of Lemma~\ref{lm:1} we have $m=\pi(\overline{\psi}(k+k'))$.) 
Since we do not have to specify  the arguments of the functions $cl_{v,n,q}[\ldots]$
the resulting recurrence relation for the cardinalities of the set of clusters is:
\begin{equation}
\label{one_binom}
 cl_{v\stackrel{\pi}{\mapsto}v',n,q} = \binom{n-m}{l-m} cl_{v',n-l+k',q-1}.
\end{equation}
 The case $l\leq k+k'$ is also covered by   \eqref{one_binom}  
since in this case the maximal entry $m$ should be equal to the length $l$ 
of the permutation $\pi$ and the corresponding binomial coefficient is equal to $1$.

 Introduce the next  family of  generating functions:
$$y_{v}(x,t):= \sum_{n,q} {cl}_{v,n,q} \frac{x^n}{n!} t^q, 
$$
one for each   vertex $\oo{v}$  of  $\mathcal{G}(\varPi)$.
Let $\oo{v}\stackrel{\pi_1}{\mapsto} \ooi{v_1},$\ldots, $\oo{v}\stackrel{\pi_d}{\mapsto} \ooi{v_d}$ 
 be the list of all edges in  $\mathcal{V}(\varPi)$ starting at $\oo{v}$.  
Denote by $l_j$ the length of the pattern $\pi_j$; denote by $k_j$  the length of  $v_j$,  and, finally,  
denote by $m_j$ the  maximal entry among  the $k_j$ last  entries of $\pi_j$. 
Formula \eqref{one_binom} implies the following reccurence relation:
\begin{equation}
\label{monot_recurtion}
 {cl}_{v,n,q} = \sum_{j=1}^d \binom{n-m_j}{l-m_j} {cl}_{v_j,n-l_j+k_j,q-1};\quad
% \text{   with the initial values   }�
{cl}_{v,n,0}=\begin{cases} 1,\; v=(1)\\ 0,\; v\neq (1).\end{cases}
\end{equation}

\begin{theorem}
\label{th:monotone_patterns}
Given a monotone collection of patterns  $\varPi$ one has that  the cluster generating function $\varPi_{cl}(x,t)=y_{(1)}(x,t)$ together with all $y_{v}(x,t),\, \oo{v}\in \mathcal{V}(\varPi)$ solve the following  system  of  linear ordinary differential equations in $x$:  
\begin{equation}
\label{monot_dif_eq}
 \frac{d^m}{d x^m}y_{v}(x,t) = t \sum_{j=1}^d 
\frac{d^{m-m_j}}{d x^{m-m_j}} \left(
\frac{x^{l_j-m_j}}{(l_j-m_j)!} \frac{d^{k_j}}{d x^{k_j}}
 y_{{v_j}}(x,t)\right). 
\end{equation}
Here $m := max\{m_j\}$  and $\oo{v}$ runs over the set $\mathcal{V}(\varPi)$ of all vertices of  $\mathcal{G}(\varPi)$.  
(Boundary conditions for each $y_v(x,t)$ can be easily determined in each particular case using  the initial terms in \eqref{monot_recurtion}.) 
\end{theorem}

 \begin{proof}
Follows  from \eqref{monot_recurtion}.   
\end{proof}

\begin{remark} 
As an immediate consequence of Theorem~\ref{th:monotone_patterns} one gets that $\varPi_{cl}(x,t)=y_{(1)}(x,t)$ satisfies a certain high order linear ordinary differential equation with polynomial coefficients (which can be obtained from the above system after  the  elimination of all $y_{v}(x,t),\; \oo{v}\neq \one$.)  
\end{remark}�

In particular,  one can get the following simplification of the Theorem~\ref{th:main} for monotone collections of patterns.
\begin{corollary}
Two monotone collections of patterns $\varPi_1$ and $\varPi_2$ are strongly c-Wilf equivalent 
if there exists a bijection $\varphi: \varPi_1 \rightarrow \varPi_2$  preserving the first two  properties as in Theorem~\ref{th:main} (i.e. preserving  lengths and linkages)
and, additionally,   preserving the maxima of the overlapping sets.
 More precisely,  for any pair  $(\pi,\pi')$ of patterns $\pi,\pi'\in\varPi_1$ of lengths $l$ and $l'$ resp.  
and an arbitrary positive integer  $k\leq min(l,l')$ the 
coincidence of the standardizations 
$\st[(\pi(l-k+1)\ldots\pi(l))] = \st[(\pi'(1)\ldots\pi'(k))]$
 implies the coincidence of the maxima of the sets:  
$$ max\{\pi(l-k+1),\ldots,\pi(l)\}= max\{\varphi(\pi)(l-k+1),\ldots,\varphi(\pi)(l)\}.$$
(Two other sets $\{\pi'(1),\ldots,\pi'(k)\}$ and $\{\varphi(\pi')(1)\ldots\varphi(\pi')(k)\}$ should coincide with the standard set 
$\{1,\ldots,k\}$ because of the monotonicity property.) 
\end{corollary}
With this in mind one can add two more  strongly c-Wilf equivalent patterns to the pair~\eqref{wilf_example} discussed in the introduction.
Namely,  
$$
(143265987) \sim_{\text{c-W}} (134265897) \sim_{\text{c-W}} (143256987) \sim_{\text{c-W}} (134256897).  
$$

\subsection{Case of a single monotone pattern}
Let us present further  simplifications for monotone collections with  a single pattern.
As it was shown in \S~\ref{section-one-pattern} we do not really need 
the graph        $\mathcal{G}(\varPi)$      for collections with one pattern.
Take  a monotone pattern $\pi$ of length $n$ and 
let $1=k_1<\ldots<k_d=k$ be the list of all different lengths of self-overlappings of $\pi$. 
Let $\pi(l)=m_1\leq \ldots \leq m_d = m$ be the list of maxima among the corresponding number of  
the last entries of $\pi$, i.e.  
 $m_j:= max\{\pi(l-k_j+1),\ldots,\pi(l)\}$.
\begin{corollary}
Given a monotone pattern $\pi$ one has the following recurrence formula for the set of its $q$-clusters of length $n$:
$$
 {cl}_{n,q} = \sum_{j=1}^d \binom{n-m_j}{l-m_j} {cl}_{n-l_j+k_j,q-1}
$$
and its  cluster generating function $y_{(1)}(x,t)=\varPi_{cl,\{\pi\}}(x,t)$  solves the linear ordinary differential equation:
\begin{equation}
\label{monot_dif_eq_one_pat}
 \frac{d^m}{d x^m}y_{(1)}(x,t) = t \sum_{j=1}^d 
\frac{d^{m-m_j}}{d x^{m-m_j}} \left(
\frac{x^{l-m_j}}{(l-m_j)!} \frac{d^{k_j}}{d x^{k_j}}
 y_{(1)}(x,t)\right),
\end{equation}
with the boundary conditions: $y_{(1)}(0,t)=0,\; y_{(1)}'(0,t) = 1,\;y_{(1)}^{\prime\prime}(0,t)=...=y_{(1)}^{(m-1)}(0,t)=0$.
\label{lm:one-pat-dif-eq}
\end{corollary}

Two particular cases covered by this corollary, namely the pattern $(12\ldots l)$ and an arbitrary pattern of the form  $(1\ldots a)$ having no nontrivial self-overlappings were considered earlier in \cite{ElN}.

Finally, using these considerations we can completely describe which permutations of length $5$ ($\pi\in S_5$) are c-Wilf equivalent. Our  main theorem~\ref{th:main} gives necessary and sufficient condition for 
c-Wilf equivalence of two patterns of length $5$. Notice that   there are two natural transformations of patterns preserving the cluster generating functions. The first one is the  reversion that rewrites a pattern backwards, i.e. it sends ($(\pi(1)\ldots\pi(n))$ to $(\pi(n)\ldots\pi(1))$). The second one  takes the complement of a pattern, i.e.  it sends $((\pi(1)\ldots\pi(n))$ to  $(n-\pi(1)+1,\ldots, n-\pi(n)+1)$.
They generate the group $\mathbb{Z}_2\times\mathbb{Z}_2$ acting on $S_5$. One can easily check that this action has $32$ orbits of which $4$ orbits with representatives 
$ 12345, 14325, 21354, 25314$ have length $2$ and the remaining $28$ have length $4$. Additionally, there are 14 orbits  whose permutations have no nontrivial selp-overlappings; 15 orbits with the only nontrivial self-overlapping of length $2$; 2 orbits with the only nontrivial self-overlapping of length $3$, and a single orbit with self-overlappings of length 2, 3 and 4, see the lists of representatives  in Proposition~\ref{pr:S5}.  

\begin{proposition}\label{pr:S5} 
Subdividing  the representatives of the orbits of the latter group action on $S_5$  into the subsets  according to the lengths of their  maximal possible self-overlappings one gets the following: 
\begin{itemize}
\item
c-Wilf equivalent orbits   having no  nontrivial overlappings are enumerated by their first and last elements:
$$
13452\sim_{c-W} 13542\sim_{c-W} 14352\sim_{c-W} 14532\sim_{c-W} 15342\sim_{c-W} 15432; $$ 
$$ 12453\sim_{c-W} 12543; \quad 12354\sim_{c-W}13254;\quad  21354\sim_{c-W} 21534; \quad 24153\sim_{c-W} 25143.
$$
\item 
among all $15$ orbits  with only 2-overlappings  having the representatives  
$$
12435, 12534, 13425, 13524, 14325, 14523, 15324, 15423, 15234, 21453, 21543, 23514, 24513, 25314, 25413
$$
 no two are  strongly c-Wilf equivalent to each other; 
\item $2$ orbits with representatives $14253, 15243$ having  $3$-overlappings are not c-Wilf equivalent; 

\item
the unique orbit with overlappings of lengths $2,3,4$ is represented by the  monotone pattern $(12345)$.
\end{itemize}
\end{proposition}

\begin{proof} Besides the application of   Theorem~\ref{th:main} one has to check by hand that the number of $3$-clusters is distinguishing all the 15 orbits  with only 2-overlappings and already the number of $2$-clusters distinguishes between the two orbits having only $3$-overlappings. 
\end{proof}�

\subsection{Examples of differential equations}
 Let us finish the  paper by presenting the (system of) linear differential equations  for our examples in Fig.\picref{graph-example}.  
One can easily check that all the collections  $\varPi_1 - \varPi_4$ shown there are monotone. 
Set $\pi_1=(1342765)$ and $\pi_2=(152364)$ for $\varPi_1$ and notice that its graph contains a single vertex. Then in the above notation one has $l_1=7,\; m=m_1=5,\; l_2=6,\; m_2=4,\; k=1$. Equation \eqref{monot_dif_eq} for the cluster generating function $y_{(1)}(x,t)$ then takes the form:
 $$y^{\rm{V}}_{(1)}=t\left(\frac{x^2}{2} y^\prime_{(1)} +\frac {d}{dx}\left(\frac{x^2}{2}y^\prime_{(1)}\right)\right).$$
 For the collection  $\varPi_2=\{(132679485)\}$ that has self-overlappings of lengths $1$ and $3$ 
one has $l_1=9,\; k_1=1,\; k=k_2=3,\; m_1=5,\; m=m_2=8$. With this data equation \eqref{monot_dif_eq_one_pat} takes the form:
$$y^{\rm{VIII}}_{(1)}=t\left(\frac{d^3}{dx^3}\left(\frac{x^3}{3!}y^\prime_{(1)}\right)+xy^{\prime\prime\prime}_{(1)}\right).$$
For the collection $\varPi_3=\{1576243, 13254\}$  denote its edges by 
$A=(\{1\},\{3\},7)$, $B=(\{1\},\{4\},5)$, $C=(\{1\},\{2,3,4\},7)$ and $D=(\{1,2,3\},\{4\},5)$.  Then one has $l_A=7, k_A=1, m_A=3$; $l_B=5, k_B=1, m_B=4$; $l_C=7, k_C=3, m_C=4$;  $l_D=5, k_D=1, m_D=4$. Thus, one gets the following system of equations:
$$
\left\{
\begin{array}{l}
 y_{(1)}^{\rm{IV}}=t\left( \frac{d}{dx} \left( \frac{x^4}{4!} y_{(1)}^{\prime}\right) +
xy_{(1)}^\prime+\frac{x^3}{3!}y_{(132)}^{\prime\prime\prime}\right), \\
 y_{(132)}^{\rm{IV}} = tx y_{(1)}^\prime.
\end{array}
\right.
$$

In this case it is easy to get a linear ordinary differential  equation satisfied by $y_{(1)}$ by dividing both sides of the first equation by $\frac{x^3}{3!}$, differentiating the result with respect to $x$  and equating the expressions for  $y_{(132)}^{\rm{IV}}$ from the first and second equations. The resulting equation has the form
$$6\left(x^3y_{(1)}^{\rm {V}}-3x^2y_{(1)}^{\rm {IV}}-tx^3\left(\frac{x^4}{4!}y'_{(1)}\right)^{\prime\prime}+3tx^2\left(\frac{x^4}{4!}y'_{(1)}\right)^{\prime}-tx^3(xy'_{(1)})^{\prime}+3tx^3y_{(1)}^\prime\right)-tx^7y_{(1)}^\prime=0.$$

Finally, for the collection $\varPi_4=\{(12534), (132465)\}$ 
denote its edges by $A=(\{1\},\{4\},5)$; $B=(\{1\},\{5\},6)$; $C=(\{1\},\{3,4,5\},5)$; $D=(\{1\},\{4,5,6\},6)$; 
$E=(\{1,2,3\},\{5\},6)$, and  $F=(\{1,2,3\}, \{4,5,6\},6)$. 
Then one has $l_A=5, k_A=1, m_A=4$; 
$l_B=6, k_B=1, m_B=5$; $l_C=5, k_C=3, m_C=5$;  $l_D=6, k_D=3, m_D=6$; $l_E=6, k_E=1, m_E=5$; $l_F=6, k_F=3, m_F=6$. 
Thus, one gets the following system of equations:
$$
\left\{
\begin{array}{l}
 y_{(1)}^{\rm{VI}}=t\left(\frac{d^2}{dx^2}\left( xy_{(1)}^{\prime}\right) +\frac{d}{dx}\left(xy_{(1)}^\prime\right)+y^{\rm{IV}}_{(132)}+y^{\rm{\prime\prime\prime}}_{(132)}\right),\\
 y_{(132)}^{\rm{VI}} = t\left( \frac{d}{dx} \left(x y_{(1)}^\prime \right)+y_{(132)}^{\prime\prime\prime}\right). 
\end{array}
\right.
$$ 

\medskip
Let us extract from this system a linear ordinary differential equation satisfied by $y_{(1)}$. To simplify our notation set $u=y_{(1)}^\prime,\; v=y_{(132)}^{\prime\prime\prime}.$ Then  we get the system: 
$$
\left\{
\begin{array}{l}
u^{\rm{V}}=t\left(\left( xu\right)^{\prime\prime} +\left(xu\right)^\prime+v^\prime+v\right),\\
v^{\prime\prime\prime} = t\left( (xu)^\prime+v\right). 
\end{array}
\right.
$$
Since the coefficients in the left-hand sides are both equal to $1$ we can eliminate $v$ by differentiating the first equation a number of times and substituting $v^{\prime\prime\prime}$ from the second equation till we can solve both equations for the remaining $v$ and its derivative of appropriate order. In our concrete case, differentiating the first equation twice and substituting $v^{\prime\prime\prime}$ from the second equation we get
$$
\left\{
\begin{array}{l}
u^{\rm{VII}}=t(( xu)^{\rm{IV}} +(xu)^{\prime\prime\prime}+t(xu)^\prime+tv+v^{\prime\prime}),\\
v^{\prime\prime\prime} = t\left( (xu)^\prime+v\right). 
\end{array}
\right.
$$
Differentiating the first equation in the latter system once and substituting  $v^{\prime\prime\prime}$ again we get
$$
\left\{
\begin{array}{l}
u^{\rm{VIII}}=t(( xu)^{\rm{V}} +(xu)^{\rm{IV}}+t(xu)^{\prime\prime}+t(xu)^\prime+ t(v^\prime+v)),\\
v^{\prime\prime\prime} = t\left( (xu)^\prime+v\right). 
\end{array}
\right.
$$
Finally, equating the expressions for $v+v^\prime$ from the first equation in the latter system and from the original equation for $u^{\rm{V}}$ we get the following equation 
$$u^{\rm{VIII}}-t(xu)^{\rm V}-tu^{\rm{V}}-t(xu)^{\rm{IV}}+t^2(xu)^{\prime\prime\prime}-t^2(xu)^{\prime\prime}+t^3(xu)^\prime-t^2(xu)^\prime=0$$
containing $u$ and its derivatives only. Substituting $u=y^\prime_{(1)}$ we obtain the required equation 
$$  y^{\rm{IX}}_{(1)}-t(xy^\prime_{(1)})^{\rm V}-ty^{\rm{VI}}_{(1)}-t(xy^\prime_{(1)})^{\rm{IV}}+t^2(xy^\prime_{(1)})^{\prime\prime\prime}-t^2(xy^\prime_{(1)})^{\prime\prime}+t^3(xy^\prime_{(1)})^\prime-t^2(xy^\prime_{(1)})^\prime=0
 $$ 
for the cluster generating function. 

\begin{remark}
Notice that since the leading terms in the left-hand sides of the system \eqref{monot_dif_eq} are always equal to $1$ the elimination process similar to the one just described will always lead to an equation satisfied by the cluster generating function $\varPi(x,t)$. On the other hand, there is no guarantee that the obtained linear ordinary equation with polynomial coefficients will have the minimal possible order  among such equations satisfied by $\varPi(x,t)$.
\end{remark}

\end{document}